\documentclass{amsart}
\usepackage{amssymb,   latexsym, ulem, amscd} 
\usepackage[all]{xy}

\newtheorem{theorem}{Theorem}[section]
\newtheorem{lemma}[theorem]{Lemma}
\newtheorem{sublemma}[theorem]{Sublemma}
\newtheorem{corollary}[theorem]{Corollary}
\newtheorem{proposition}[theorem]{Proposition}

\theoremstyle{definition}

\theoremstyle{conjecture}
\newtheorem{conjecture}[theorem]{Conjecture}

\theoremstyle{remark}
\newtheorem{remark}[theorem]{Remark}

\newcommand{\lem}[2]{\begin{lemma}\label{#1} #2\end{lemma}}

\newcommand{\prop}[2]{\begin{proposition}\label{#1}#2\end{proposition}}
\newcommand{\propr}[3]{\begin{proposition}[\mbox{#2}]\label{#1}#3\end{proposition}}
\newcommand{\thm}[2]{\begin{theorem}\label{#1}#2\end{theorem}}

\newcommand{\cor}[2]{\begin{corollary}\label{#1}#2\end{corollary}}

\renewcommand{\labelenumi}{\arabic{enumi})}
\newcommand{\al}{\alpha}
\newcommand{\be}{\beta}
\newcommand{\ga}{\gamma}

\newcommand{\ze}{\zeta}
\newcommand{\et}{\eta}

\newcommand{\io}{\iota}
\newcommand{\ka}{\kappa}

\newcommand{\Si}{{\Sigma}}

\newcommand{\ph}{\varphi}

\newcommand{\AR}[2]{$$\begin{array}{#1}#2\end{array}$$}

\newcommand{\ARr}[3]{$$\begin{array}{#1}#2\end{array}\lnr{#3}$$}

\newcommand{\naka}[1]{\begin{center}#1\end{center}}

\newcommand{\cass}[1]{\begin{cases}#1\end{cases}}
\newcommand{\ak}{\quad}

\newcommand{\AL}[1]{\begin{align*}#1\end{align*}}

\newcommand{\xycr}[3]{$$\xymatrix@C=#1pt@R=#2pt{#3}$$} 
\newcommand{\xycrr}[4]{$$\xymatrix@C=#1pt@R=#2pt{#3}\lnr{#4}$$} 
\newcommand{\xyr}[2]{$$\xymatrix@C=20pt@R=#1pt{#2}$$} 
\newcommand{\xyrr}[3]{$$\xymatrix@C=20pt@R=#1pt{#2}\lnr{#3}$$} 
\newcommand{\lxyr}[2]{$\xymatrix@C=20pt@R=#1pt{#2}$} 

\renewcommand{\Xy}[2]{$$\xymatrix@=#1pt{#2}$$}
\newcommand{\Xyr}[3]{$$\vcenter{\xymatrix@=#1pt{#2}}\lnr{#3}$$}
\newcommand{\xR}[2]{$\xymatrix@R#1pt{#2}$}
\newcommand{\xyR}[2]{$$\xymatrix@R#1pt{#2}$$}
\newcommand{\xyRr}[3]{$$\vcenter{\xymatrix@R#1pt{#2}}\lnr{#3}$$}

\newcommand{\Z}{\mathbb Z}
\newcommand{\Zt}{\mathfrak Z}
\newcommand{\N}{\mathbb N}

\newcommand{\Zp}{{\mathbb Z}_{(p)}}
\def\o+{\oplus}
\newcommand{\Op}{\bigoplus}

\newcommand{\ox}{\otimes}

\newcommand{\e}{{\rm Ext}}
\newcommand{\lnb}{\refstepcounter{theorem}\leqno(\thetheorem)}
\newcommand{\lnr}[1]{\lnb\label{#1}}

\newcommand{\nr}{\refstepcounter{theorem}\thetheorem}
\newcommand{\kko}[1]{(\ref{#1})}

\newcommand{\mbx}[1]{\quad\mbox{#1}\quad}
\newcommand{\cf}{{\it cf.}\ }
\newcommand{\sus}{\subset }
\newcommand{\sm}{\wedge }

\newcommand{\Ker}{{\rm Ker}\ }

\newcommand{\im}{{\rm Im}\ }

\newcommand{\cg}{\equiv}
\newcommand{\lrk}[1]{\left\langle #1\right\rangle}
\newcommand{\qand}{\mbx{and}}

\newcommand{\Lt}{\left}
\newcommand{\Rt}{\right}
\newcommand{\LR}[1]{\Lt(#1\Rt)}
\newcommand{\cln}{\colon}
\renewcommand{\O}[1]{\overline{#1}}

\newcommand{\mx}[1]{\begin{matrix}#1\end{matrix}}

\renewcommand{\b}[2]{\ensuremath{b_{#1,#2}}}
\newcommand{\h}[2]{\ensuremath{h_{#1,#2}}}
\newcommand{\el}[2]{\ensuremath{\ell_{#1,#2}}}

\newcommand{\pic}[1]{\mbox{{\rm Pic}}(\cL_{#1})}

\newcommand{\pc}[1]{\mbox{{\rm Pic}}^0(\cL_{#1})}
\newcommand{\thick}[1]{\mbox{{\rm thick}}\lrk{#1}}

\newcommand{\cI}{{\mathcal I}}
\newcommand{\cK}[1]{{S_{(#1)}}}
\newcommand{\cL}{{\mathcal L}}

\newcommand{\cP}{{\mathcal P}}
\newcommand{\cS}{{\mathcal S}}
\newcommand{\cT}{{\mathcal T}}
\newcommand{\xar}{\xrightarrow}
\newcommand{\xlar}{\xleftarrow}

\newcommand{\p}{\noindent {\it Proof.}  }
\newcommand{\q}{\hfill \qed \\[-.05in] }

\newcommand{\LK}{L_{K(n)}}

\newcommand{\holim}{\mathop{\rm holim}\limits}

\newcommand{\rank}{\mbox{\rm rank}\, }

\newcommand{\MJ}[1]{MJ^{(#1)}}
\newcommand{\jo}[1]{^{(#1)}}

\newcommand{\llb}[1]{{[\![} #1{]\!]}}

\begin{document}

\title[Generalized Moore spectra and Hopkins' Picard groups]{ Generalized Moore spectra and Hopkins' Picard groups for a smaller chromatic level} 
\author{Ryo Kato}
\address{Faculty of Fundamental Science, National Institute of Technology, Niihama College, Niihama, 792-8580, Japan}
\email{ryo\_kato\_1128@yahoo.co.jp}
\author{You-na Kawamoto}
\address{Shin'ei-cho, Kagoshima, 890-0072, Japan}
\author{Hiroki Okajima}
\address{Department of Mathematics, Faculty of Science, Kochi University, Kochi,
780-8520, Japan}
\email{gg1122cc@gmail.com}
\author{ Katsumi Shimomura}
\address{Department of Mathematics, Faculty of Science, Kochi University, Kochi,
780-8520, Japan}
\email{katsumi@kochi-u.ac.jp}
\begin{abstract}
Let $\cL_n$ for a positive integer $n$ denote the stable homotopy category of $v_n^{-1}BP$-local spectra at a prime number $p$.
Then, M.~Hopkins defines the Picard group of $\cL_n$ as a collection of isomorphism classes of invertible spectra,
whose exotic summand $\pc n$ is studied by several authors.
In this paper, we study the summand for $n$ with $n^2\le 2p+2$.
For $n^2\le 2p-2$, it consists of invertible spectra whose $K(n)$-localization is the $K(n)$-local sphere.
In particular, $X$ is an exotic invertible spectrum of $\cL_n$ if and only if $X\sm MJ$ is isomorphic to a $v_n^{-1}BP$-localization  of the generalized Moore spectrum $MJ$ for an invarinat regular ideal $J$ of length $n$.
For $n$ with $2p-2<n^2\le 2p+2$, we consider the cases for $(p,n)=(5,3)$ and $(7,4)$.
In these cases, we characterize them by the Smith-Toda spectra $V(n-1)$.
For this sake, we show that $L_3V(2)$ at the prime five and $L_4V(3)$ at the prime seven are ring spectra.
\end{abstract}
\maketitle

\section{Introduction}

Let $\mathcal S_p$ be the stable homotopy category of $p$-local spectra for an odd prime number $p$.
Consider the Brown-Peterson spectrum $BP$ characterized by the homotopy groups $\pi_*(BP)=BP_*=\Zp[v_1,v_2,\dots]$ over generators $v_k$
with degree $|v_k|=2(p^k-1)$.
We work in the stable homotopy category $\cL_n$ for $n\ge 0$ consisting of $v_n^{-1}BP$-local spectra.
A spectrum $X\in\cL_n$ is called {\it invertible} if there is a spectrum $Y$ such that $X\sm Y\simeq L_nS^0$.
Here, $L_n\cln \cS_p\to \cL_n$ denotes the Bousfield localization functor.
Mike Hopkins introduced the Picard group $\pic n$ consisting of the isomorphism classes 
of invertible spectra (\cf \cite{st}).
Mark Hovey and Hal Sadofsky \cite{hs} showed that $\pic n$ is an abelian group with the decomposition 
$$
\pic n\cong \Z\o+ \pc n, \lnr{dec}
$$
and 
$$
\pc n=0 \ \mbx{if $n^2+n\le q$} \lnr{pic0}
$$
for the integer
$$
q=2p-2.
$$

For each $n\ge 1$, 
consider an invariant ideal
$$
J=(p^{e_0},v_1^{e_1},\dots,v_{n-1}^{e_{n-1}}) \sus BP_* \lnr{inv}
$$
for positive integers $e_i$ (see \kko{tsu}).
We call a spectrum $MJ$  with $BP_*(MJ)=BP_*/J$ a {\it type $n$ generalized Moore spectrum}.

\medskip

\noindent
(\nr) \label{hs-d} (Hopkins and Smith \cite{hsm}) {\it For each invariant ideal $J$ of the form \kko{inv}, there exists a type $n$ generalized Moore spectrum $MJ'$ for an invariant ideal $J'$ of the form \kko{inv} with $J'\sus J$.\\}

Let $MJ$ for an invariant  ideal $J$ in \kko{inv} be a type $n$ generalized Moore spectrum, and put
$$
S_J=\{[X]\mid X\in \thick{L_nS^0},\ X\sm MJ\simeq L_nMJ\}, 
$$
where $[X]$ denotes the isomorphism class of $X$.

\thm{SJ}{
$S_J$ is a subgroup of $\pc n$.
}

Let $E(n)$ be the $n$-th Johnson-Wilson spectrum. Then, the category $\cL_n$ also consists of $E(n)$-local spectra.
The spectrum gives rise to
 the Hopf algebroid
\AL{
(E(n)_*&,E(n)_*(E(n)))\\
&=(\Zp[v_1,v_2,\dots,v_n,v_n^{-1}], E(n)_*\ox_{BP_*}BP_*(BP)\ox_{BP_*}E(n)_*)
}
induced from the Hopf algebroid $BP_*(BP)=BP_*[t_1,t_2,\dots]$ over $t_k$ with $|t_k|=2(p^k-1)$.
We notice that 

\medskip

\noindent
(\nr) \label{ks} (\cite[Th.~1.1]{ks}) {\it $X\in \pc n$ if and only if $E(n)_*(X)\cong E(n)_*$ as an $E(n)_*(E(n))$-comodule.}

\medskip

We put
$$
H^{s,t}M=\e_{E(n)_*(E(n))}^{s,t}(E(n)_*,M) \mbx{for an $E(n)_*(E(n))$-comodule $M$.} \lnr{homology}
$$
Then, we have
the $E(n)$-based Adams spectral sequence
$$
E_2^{s,t}(X)=H^{s,t}E(n)_*(X)\Longrightarrow \pi_{t-s}(L_nX). \lnr{EnASS}
$$
The isomorphism in \kko{ks} induces an isomorphism
$$
E_2^{s,t}(X)\cong E_2^{s,t}(S^0)\mbx{ for $X\in \pc n$.}\lnr{E2}
$$

From now on,  we assume that the integer $n$ satisfies
$$
n\ge 3 \qand n^2+n\le 2q. \lnr{cond}
$$
We notice that $n<p-1$ under \kko{cond}, and that by \cite[(10.10)]{r-loc},
$$
E_2^{kq+1,kq}(S^0)=0 \mbx{for $k>1$.} \lnr{E20}
$$
In this case, we have a monomorphism $\ph\cln \pc n\to  E_2^{q+1,q}(S^0)$ by \cite[Th.~1.2]{ks} defined by $\ph(X)=w$ for $w$ in the differential $$
d_{q+1}(1_X)=w1_X\in E_2^{q+1, q}(X) ,\lnr{d1}
$$ 
where $1_X$ is the generator of $E_2^{0,0}(X)\cong E_2^{0,0}(S^0)\cong \Z_{(p)}$.
The monomorphism is actually an isomorphism:

\medskip

\noindent
(\nr)\label{P=E} (\cite[Cor.~1.9]s) {\it Under the condition \kko{cond}, $\pc n\cong E_2^{q+1,q}(S^0)$.}

\medskip

Let $MJ$ be a type $n$ generalized Moore spectrum for an invariant ideal $J$ of \kko{inv}, and  
$$
i_J\cln S^0\to MJ, \lnr{iJ}
$$ 
the inclusion  to the bottom cell.
Consider the induced homomorphism 
$$
(i_J)_*\cln E_2^{q+1,q}(S^0)\to E_2^{q+1,q}(MJ). \lnr{iJ*}
$$
Then, Theorem \ref{SJ} and \kko{P=E} imply the following: 

\prop0{ Suppose that an integer $n$ satisfies \kko{cond}. Then, 
$$
\Ker (i_J)_*\cong S_J
$$ 
for $(i_J)_*$ in \kko{iJ*} if $L_nMJ$ is a ring spectrum.
}

This together with \kko{P=E} implies the following theorem:

\thm{KerJ}{Let $n$ be an integer satifying \kko{cond}, and suppose that there exists a type $n$ generalized Moore spectrum $MJ$ for an ideal $J$  
in \kko{inv}. If $L_nMJ$ is a ring spectrum with $E_2^{q+1,q}(MJ)=0$, then
$\pc n=E_2^{q+1,q}(S^0)=\Ker (i_J)_*=S_J$.
}

In this paper, we call $M$ a {\it ring spectrum} if there exist maps $\mu_M\cln M\sm M\to M$ and $i_M\cln S^0\to M$ such that the composite $M=S^0\sm M\xar{i_M\sm M}M\sm M\xar{\mu_M}M$ is homotopic to the identity. 
For a spectrum $MJ'$ in \kko{hs-d}, Devinatz further showed

\medskip

\noindent
(\nr) \label{d:ring} (Devinatz \cite d) {\it We may take $MJ'$ in \kko{hs-d} to be a ring spectrum.}

\medskip

The celebrated theorems \kko{hs-d} and \kko{d:ring} enable us to consider an inverse system of type $n$ generalized Moore ring spectra
$$
\MJ 1\xlar{\pi^1}\MJ2\xlar{\pi^2} \cdots \xlar{\pi^{k-1}}\MJ k\xlar{\pi^{k}}\MJ{k+1}\xlar{\pi^{k+1}}\cdots,
$$
in which $\ J^{(k)}\supset J^{(k+1)}$ are invariant ideals such that $\bigcap_{k\ge 1}J^{(k)}=0$, and
 $BP_*(\pi^{k})\cln$ $ BP_*(\MJ{k+1})\to BP_*(\MJ{k})$ are the canonical projections.
Furthermore, we assume that the Spanier-Whitehead dual $D(\MJ k)$ of $\MJ k$ is isomorphic to $\Si^a \MJ k$ for some $a\in\Z$.
We fix such an inverse system, and denote the set of the invariant ideals by
$$
\cI_n=\{J^{(k)}\sus BP_*\mid \mbox{$J^{(k)}$ is the invariant ideal in the above system}\}. 
\lnr{mn}
$$

Let $L_E\colon \cS \to \cS$ denote the Bousfield localization functor with respect to a spectrum $E$. We notice that $L_n=L_{v_n^{-1}BP}$. Furthermore, $L_{F(n)}$ denotes $L_F$ for a type $n$ finite spectrum $F$, which is well defined since $L_F=L_{F'}$ for type $n$ finite spectra $F$ and $F'$.
$$
L_{F(n)}X\simeq \holim _{J\in \cI_n}X\sm MJ\leqno{(\nr) \label{Fnloc} \mbox{(Hovey \cite[Th.~2.1]{h})}}
$$

Consider a spectrum $v_n^{-1}BP\sm MJ$ for $J$ of \kko{inv}.
Then, the Bousfield class $\lrk{v_n^{-1}BP\sm MJ}$ of $v_n^{-1}BP\sm MJ$ equals the Bousfield class $\lrk {K(n)}$ of the $n$-th Morava $K$-theory $K(n)$. 
Since  $L_{F(n)\sm E}X=L_{F(n)}L_EX$ by \cite[Cor.~2.2]h,
we obtain $L_{F(n)}L_nX\simeq L_{K(n)}X$.
In particular, taking $X=L_nS^0$ in \kko{Fnloc}, we have
$$
L_{K(n)}S^0\simeq \holim_{J\in\cI_n}L_nMJ.
$$
Consider the kernel of the homomorphism $\pc n\to \ka_n \sus \pic{K(n)}$ (\cf \cite[Cor.~2.5]{hs}) induced from the localization $L_{K(n)}\cln \cL_n\to \cL_{K(n)}$:
$$
\cK n=\{[X]\in \pc n\mid L_{K(n)}S^0\simeq L_{K(n)}X\}. \lnr{kerK}
$$
\prop{cap}{$\bigcap_{J\in \cI_n}S_J=\cK n$.}

The decomposition \kko{dec} implies that for any invertible spectrum $X$ in $\cL_n$, there exists an integer $s$ such that $\Si^s X$ represents an element of $\pc n$. Moreover,
Morava's structure theorem implies that $E_2^{q+1,q}(MJ)=0$ if $n^2\le q$ and $n\ge 3$, and so
Theorem \ref{KerJ} and Proposition \ref{cap} as well as \kko{pic0} imply the following:

\cor{C-E2}{Let $p$ and $n$ be the integers in \kko{cond}. If $n^2+n\le q$, then $\pc n=0$.
If $n^2\le q<n^2+n$,  then $\pc n=S_J=\cK n$. 
In other words, the homomorphism $\pc n\to \pic {K(n)}$ induced from $L_{K(n)}$ is the zero homomorphism if $n^2\le q$.
Furthermore, $X$ is invertible in $\cL_n$  if and only if  $X\in \thick{L_nS^0}$ and $X\sm MJ\simeq \Si^sL_nMJ$ for an integer $s$ and an ideal $J$ of length $n$ in \kko{inv}.
}

Now suppose that $q<n^2$.
In this case, we have little knowledge about the homomorphisms $(i_J)_*\cln E_2^{q+1,q}(S^0)\to E_2^{q+1,q}(MJ)$. 
We notice that there is no pair $(p,n)$ satisfying $n^2-2=q$ or $n^2-3=q$.
So we consider the cases 
\begin{enumerate}
\item
$(p,n)=(5,3)$, under which $n^2-1=q$, and
\item
$(p,n)=(7,4)$, under which $n^2-4=q$.
\end{enumerate}

Note that in these cases, the Smith-Toda spectrum $V(n-1)=MI_n$ exists but it is not a ring spectrum (\cf \cite{gb}).
Here, $I_n=(p,v_1,\dots, v_{n-1})$ is the invariant prime ideal of $BP_*$.
In this paper, we show the following:

\thm{ring}{For $(p,n)=(5,3)$ or $(7,4)$, $L_nV(n-1)$ is a ring spectrum.}

\thm{small}{For $(p,n)=(5,3)$ or $(7,4)$, $E_2^{q+1,q}(V(n-1))=0$.}

Theorems \ref{KerJ}, \ref{ring} and \ref{small} imply

\cor{C1-V2}{For $n=3,4$, we have
$$
\pc 3=\cass{S_{I_3}&p=5\\
0&p\ge 7} \qand \pc 4=\cass{S_{I_4}&p=7\\
0&p\ge 11.}
$$
}

We notice that $\pc 3$ at the prime five is isomorphic to $S_{J_k}$ for $J_k=(5,v_1,v_2^k)$ with some $k>1$ (at least $p^2$),  but the result may be less interesting and omit here. 

In the next section, we show Theorem \ref{SJ} and Propositions \ref0 and \ref{cap}.
We verify Theorem \ref{ring} in section three.
Section four is devoted to showing Theorem \ref{small}.

 \section{The group $S_J$}

Consider an ideal 
$$
J=(p^{e_0}, v_1^{s_1p^{e_1}},\dots , v_{n-1}^{s_{n-1}p^{e_{n-1}}}) \lnr{ideal}
$$
of $BP_*$, where $e_i\ge 0$, $s_i\ge 1$ and $p\nmid s_i$.
Then, we notice the following:

\medskip

\noindent
(\nr) \label{tsu}(\cite[Th.~1.5]{tsuk}) 
\ak {\it $J$ is an invariant ideal of $BP_*$ if and only if $e_0-1\le e_1$ and $s_i\le p^{e_{i+1}-e_i-e_0+1}$ for $1\le i <n$.
}

\medskip

Let $\thick{L_nS^0}$ denote the thick subcategory of $\cL_n$ generated by $L_nS^0$.
Since $\cL_n$ is a monogenic stable homotopy category, we see the following:

\medskip

\noindent
(\nr) \label{DX} (\cf \cite[Th.~2.1.3]{hps}) {\it If $X, Y\in \thick{L_nS^0}$, then so are $D(X)$ and $X\sm Y$.
Here, $D(X)$ denotes the Spanier-Whitehead dual of $X$ in $\cL_n$.}

\lem{SJ1}{Suppose that $C\in\thick{L_nS^0}$ satisfies $C\sm MJ=0$ for a type $n$ generalized Moore spectrum $MJ$ with $J$ in $\kko{ideal}$. Then, $C=0$.}

\p
Let $J_k=(p^{e_0},v_1^{e_1},\dots,v_{k-1}^{e_{k-1}})$ be an invariant ideal of $E(n)_*$ for each $0\le k\le n$ $(J_0=(0), J_n=J)$.
Suppose that 
$C\sm MJ_{k+1}= 0$.
Then, the cofiber sequence $MJ_k\xar{v_k^{e_k}}MJ_k\to MJ_{k+1}$ gives rise to an isomorphism $C\sm MJ_k\xar[\simeq]{C\sm v_k^{e_k}}C\sm MJ_k$, which implies $E(n)_*(C\sm MJ_k)=v_k^{-1}E(n)_*(C\sm MJ_k)$.
Since $C\sm MJ_k\in \thick{L_nS^0}$, $E(n)_*(C\sm MJ_k)$ is a finitely generated $E(n)_*$-module,
and hence
 $C\sm MJ_k=0$.
Inductively, we deduce $C=0$.
\q

\propr{ev}{\cf \cite[Lemma A.2.6]{hps}}{
Let $X\in\thick{L_nS^0}$ and $ev\cln D(X)\sm X\to S^0$ denote  the evaluation map.
Then, $ev\sm D(X)\cln D(X)\sm X\sm D(X)\to D(X)$ is a retraction.
}

\p
Consider the cofiber sequence
$$
D(X)\sm X\xar{ev} L_nS^0 \xar{c} C \lnr{cof:ev}
$$
of the evaluation map $ev$.
It suffices to show the map $c\sm D(X)$ trivial.

The evaluation map $ev$ defines a homomorphism
 $$
ev_W\cln [D(X),W\sm D(X)]_*\to [D(X)\sm X, W]_*
$$ 
by $ev_W(f)=(W\sm ev)(f\sm X)$. Consider the full subcategory
$$
\cT_X=\{W\in\cL_n\mid \mbox{$ev_W$ is an isomorphism}\}
$$
of $\cL_n$. Then, it is easy to see $\cT_X$ thick, and $L_nS^0\in \cT_X$.
It follows that $\thick{L_nS^0}\sus \cT_X$.
By \kko{DX} and the cofiber sequence \kko{cof:ev}, we see $C\in \thick{L_nS^0}$, and so $C\in \cT_X$.
Therefore, we have an isomorphism
$$
ev_C\cln [D(X),C\sm D(X)]_*\to [D(X)\sm X, C]_*.
$$
Since
$$
ev_C(c\sm D(X))=(C\sm ev)(c\sm D(X)\sm X)=c\circ ev=0,
$$
we obtain $c\sm D(X)=0$ as desired.
\q

\medskip

\noindent
{\it Proof of Theorem \ref{SJ}.} 
For $[X],[Y]\in S_J$, $X\sm Y\sm MJ\simeq X\sm MJ\simeq L_nMJ$, and so $[X\sm Y]\in S_J$.

We note that $D(MJ)=\Si^aMJ$ for some integer $a$. Then, 
\AL{
D(X)\sm MJ& \simeq D(X\sm D(MJ))\simeq D(\Si^aX\sm MJ)\\
&\simeq L_nD(\Si^a MJ)\simeq L_nD(D(MJ))\simeq L_nMJ.
}
It follows that $[D(X)]\in S_J$. 

For $[X]\in S_J$, we show that $X$ is invertible.
Consider the cofiber sequence \kko{cof:ev}.
Proposition \ref{ev} gives rise to the decomposition
$$
D(X)\sm X\sm D(X)\sm MJ\simeq (D(X)\sm MJ)\vee \Si^{-1}(C\sm D(X)\sm MJ),
$$
which induces an isomorphism 
$$
L_nMJ\simeq L_nMJ\vee \Si^{-1}(L_nC\sm MJ)
$$
by the above observation.
Since $E(n)_*(MJ)$ is a monogenic $E(n)_*$-module, the summand $E(n)_*(C\sm MJ)$ is zero.
Hence $L_nC\sm MJ=0$.
Note that $C\in \thick{L_nS^0}$ by \kko{DX}. Thus, Lemma \ref{SJ1} shows $C=L_nC=0$, and $X$ is  invertible in $\cL_n$ with $X^{-1}=D(X)$.

Furthermore, the inclution $i_J$ in \kko{iJ} induces the canonical projection $E(n)_*(i_J)\cln$ $ E(n)_*(X)\to E(n)_*(X)/J$ of $E(n)_*(E(n))$-comodules, and so we see $X$ exotic. 
\q

\noindent
{\it Proof of Proposition \ref{0}.} 
Let $\ph\cln \pc n\to E_2^{q+1,q}(S^0)$ denote the isomorphism in \kko{P=E}, and
consider the composite $\ph'\cln S_J\sus \pc n\xar[\cong]{\ph} E_2^{q+1,q}(S^0)$ for the inclusion in Theorem \ref{SJ}.
Then, for $[X]\in S_J$, $\ph'([X])=w$ for $w$ in $d_{q+1}(1_X)=w1_X$ by \kko{d1}. 
Note that the isomorphism $\et_{J}^X\cln X\sm MJ\simeq L_nMJ$ showing $[X]\in S_J$ induces an isomorphism of $E_2$-terms in the commutative diagram
\xyR{10}{
E_2^{*,*}(X)\ar[r]^-{(i_J)_*}\ar[d]_-{\et^X_*}^\cong&E_2^{*,*}(X\sm MJ)\ar[d]_\cong^-{(\et_J^X)_*}\\
E_2^{*,*}(S^0)\ar[r]^-{(i_J)_*}&E_2^{*,*}(MJ),
} 
where $\et^X_*$ denotes the isomorphism in \kko{E2}.
We also have the generator $1_{X\sm MJ}\in E(n)_*(X\sm MJ)\stackrel{(\et_J^X)_*}\cong E(n)_*(MJ)$, which gives rise to the permanent cycle $1_{X\sm MJ}=(\et_{J}^X)_*^{-1}(1_{MJ})\in E_2^{0,0}(X\sm MJ)$ 
for the generator $1_{MJ}\in E_2^{0,0}(MJ)$.  
Since $(i_J)_*(1_X)=1_{X\sm MJ}\in E_2^{0,0}(X\sm MJ)$, the naturality of the differential shows
\AL{
(i_J)_*(w)&=(\et_J^X)_*(i_J)_*(\et^X_*)^{-1}(w)=(\et_J^X)_*(i_J)_*(w1_X)\\
&=(\et_J^X)_*(i_J)_*(d_{q+1}(1_X))=(\et_J^X)_*d_{q+1}((i_J)_*(1_X))\\
&=(\et_J^X)_*d_{q+1}(1_{X\sm MJ})=d_{q+1}(1_{MJ})=0.
}
Therefore, the monomorphism $\ph'$ reduces to $\ph'\cln S_J\to \Ker (i_J)_*$.

For any $w\in \Ker (i_J)_*$, let $X_w$ denote an inverse spectrum such that $[X_w]\in \pc n$ and $d_{q+1}(1_{X_w})=w1_{X_w}\in E_2^{q+1,q}(X_w)$.
Send this relation under the homomorphism $(i_J)_*$, and we obtain
$$
d_{q+1}(1_{X_w\sm MJ})=(i_J)_*(d_{q+1}(1_{X_w}))=(i_J)_*(w1_{X_w})=0\in E_2^{q+1,q}(X_w\sm MJ).
$$
Thus,
$1_{X_w\sm MJ}\in E_2^{0,0}(X_w\sm MJ)$ is a permanent cycle and detects a map $i_J^X\cln S^0\to X_w\sm MJ$.
Since $L_nMJ$ is a ring spectrum, the map $i_J^X$ extends to the isomorphism $L_nMJ\simeq
X_w\sm MJ$.
Thus, $[X_w]\in S_J$, and $\ph'$ is the desired isomorphism.
\q

\noindent
{\it Proof of Proposition \ref{cap}.} 
Suppose $X\in \bigcap_{J\in\cI_n}S_J$. Then, 
 it is shown in \cite[Prop.~E]{s-p} that if $\pi_0(L_nMJ)$ is finite for each $J\in\cI_n$, then $\LK X\simeq \LK S^0$.
In our case, $H^{rq,rq}K(n)_*$ is finite and equals zero if $rq>n^2$ by \cite[(1.8)~Cor., (1.9)~Th.]{r} (see Lemma \ref{coh}). It follows that $E_2^{rq,rq}(MJ)$ is finite, and so is $\pi_0(L_nMJ)$.
Therefore, $X\in \cK n$.
The converse is trivial. \q

\section{Ring structures on the Smith-Toda spectra $L_{n}V(n-1)$}

We begin with the definition of the Smith-Toda spectra $V(k)=MI_{k+1}$ for the pairs $(p,k)$ of a prime number $p$ and a non-negative integer $k$ with  $2k< p$ and $k\le 3$. 
Here, $I_{k+1}=(p,v_1,\dots,v_k)$ denotes the invariant ideal of $BP_*$. 
The spectra $V(k)$ are defined by the cofiber sequences
\ARr {cl}{
S^0 \xar{p} S^0 \xar{i} V(0) \xar j S^1 &\mbox{for }\ p\ge 2,\\
\Si^qV(0) \xar{\al} V(0) \xar {i_1}V(1) \xar{j_1} \Si^{q+1}V(0)&\mbox{for }\  p\ge 3,\\
\Si^{q_2}V(1)\xar{\be} V(1) \xar{i_2} V(2) \xar{j_2} \Si^{q_2+1}V(1)&\mbox{for \  $p\ge 5$, \ and}\\
\Si^{q_3}V(2)\xar{\ga} V(2) \xar{i_3} V(3) \xar{j_3} \Si^{q_3+1}V(2)&\mbox{for }\  p\ge 7,
}{cofs}
 for
$$
q_k=2(p^k-1)=|v_k| \lnr{qn}
$$ ($q_1=q$).
Here, $\al\in [V(0), V(0)]_{q}$, $\be\in [V(1), V(1)]_{q_2}$ and $\ga\in [V(2), V(2)]_{q_3}$ are the well known $v_1$-, $v_2$-, and $v_3$-periodic maps due to Adams, Smith and Toda, respectively.
In particular, we have a cell decomposition
$$
V(2)=\LR{(S^0\cup_pe^1)\cup_{\al}(e^{q+1}\cup_p e^{q+2})}\cup_{\be}\Si^{q_2+1}\LR{(e^0\cup_pe^1)\cup_{\al}(e^{q+1}\cup_p e^{q+2})}
$$
and $V(3)=V(2)\cup_\ga \Si^{q_3}CV(2)$.
Put 
$$
\mathfrak Z(k)=\{i\in\Z\mid \mbox{$e^i$ is a cell of $V(k)$}\}. 
$$
As stated in \cite[p.~59]{toda},  if $\pi_{i-1}(V(k))=0$ for $i=s+a$ with  $a\in \mathfrak Z(k)$, then $V(k)\jo{s-1}\sm V(k)\to V(k)$ extends to $V(k)\jo s\sm V(k)\to V(k)$.  
Here, $W^{(i)}$ denotes the $i$-skeleton of $W$.
Consider the Adams-Novikov spectral sequence
$$
{^n\!E}_2^{s,t}=\e_{BP_*(BP)}^{s,t}(BP_*, BP_*(V(n-1))\implies \pi_{t-s}(V(n-1)).
$$
Let $\cP$ denote the dual of the Hopf algebra generated by the reduced power operations, isomorphic to
$\Z/p[\xi_1,\xi_2,\dots]\cong \Z/p[t_1,t_2,\dots]\cong BP_*(BP)/(p,v_1,v_2,\dots)$.
The isomorphism $BP_*(BP)/I_n\cong \cP$ up to dimension $<q_n$ induces another one
$$
{^n\!E}_2^{s,t}\cong \e_{\cP}^{s,t}(\Z/p, \Z/p)
$$
for $t-s<q_n$.
Toda \cite{toda} showed the following: 
$$
\rank \e_{\cP}^{s,t}(\Z/p, \Z/p)\le \rank \LR{P(\b kl)\ox H^{*,*}(U(L))}^{s,t}. \leqno{(\nr)\label{trank} \mbox{(\cite[Lemma 2.2]{toda})}}
$$
Here,  the module $H^{s,t}(U(L))$ is also determined  in \cite[p.55]{toda}   for $t-s\le (p^3+3p^2+2p+1)q-4$. 
In particular, for $t-s\le 2q_3+2q_2+2q+7=(2p^2+4p+6)q+7$ ($=1591$ if $p=7$),  $H^{s,t}(U(L))$ is additively generated by the elements in the table:
{\footnotesize
\AR{|c|c|c|c|c|c|}{\hline
1&h_0&h_1&g_0&k_0&k_0h_0\\
\hline
(0,0)&(1,1)&(1,p)&(2,p+2)&(2,2p+1)&(3,2p+2)\\
\hline\hline
h_2&h_2h_0&g_1&l_1&l_2&l_1h_1\\
\hline
(1,p^2)&(2,p^2+1)&\!\!(2,p^2+2p)\!\!&\!\!(3,p^2+2p+3)\!\!&\!\!(3,p^2+3p+1)\!\!&\!\!(4,p^2+3p+3)\!\!\\
\hline\hline
k_1&l_3&k_1h_1&l_1h_2&m_1&m_1h_0\\ 
\hline
\!\!\!(2,2p^2+p)\!\!\!&\!\!\!(3,2p^2+p+2)\!\!\!&\!\!\!(3,2p^2+2p)\!\!\!&\!\!\!(4,2p^2+2p+3)\!\!\!&\!\!\!(4,2p^2+4p+2)\!\!\!&\!\!\!(5,2p^2+4p+3)\!\!\!\\
\hline
\multicolumn{6}{c}{\mbox{\tiny }}\\
\multicolumn{6}{c}{\mbox{\normalsize Table \nr \label{table}
}}}
}%
In the table, the pair of integers under each element shows  the dimension of it and  the degree of it divided by $q$.

In the following, we study homotopy groups $\pi_*(V(n-1))$ based on the fact given by \kko{trank}:

\medskip

\noindent
(\nr) \label{piV} The homotopy group $\pi_{t-s}(V(n-1))$ is a subquotient of $(P(\b kl)\ox H^{*,*}(U(L)))^{s,t}$.

\subsection{The case for $(p,n)=(5,3)$:}

The set of dimensions of cells of $V(2)$ is
$$
\Zt(2)=\{0,1,9,10,49,50,58,59
\}.
$$
By \cite[Th.~4.4]{toda}, there exists the pairing 
$$
V(1)\sm V(2)\to V(2). \lnr{pairing}
$$
We notice that this follows from the fact $\pi_{i-1}(V(2))=0$ for $i=s+a$ with $s\in \{0,1,9,10\}$ and $a\in \Zt(2)$.
So we consider the homotopy groups $\pi_{i-1}(V(2))$ for $i=s+a$ with $s\in \{49,50,58,59\}$ and $a\in \Zt(2)$.

By \kko{trank} together with Table \ref{table},
the homotopy groups of degrees $\le 121$ are subquotients generated by the following elements:
\AR{|c||c|c|c|c|c|c|c|c|}{
\hline
\deg&0& 7& 38& 39&45&54&76&77\\ 
\hline
&1&\h10&\b10&\h11&\h10\b10&g_0&\b10^2&\h11\b10\\
\hline
\deg&83&86&92&93&114 &115&121\\ 
\cline{1-8}
&\h10\b10^2&k_0&g_0\b10&k_0h_0&\b10^3&\h11\b10^2&\h10\b10^3\\
\cline{1-8}
}
This together with \kko{piV} shows that $\pi_{i-1}(V(2))=0$ for $i=s+a$ with $s\in \{0,1,9,10, 49, 50\}$ and $a\in \Zt(2)$, and
$\pi_{115}(V(2))$ is a subquotient of  $\Z/5\{\h11\b10^2\}$. 
Thus, by the next lemma, the pairing \kko{pairing} extends to $V(2)\sm V(2)\to L_3V(2)$ as desired.

\lem{htV(2)}{The homotopy groups $\pi_{115}(L_3V(2))$, $\pi_{116}(L_3V(2))$ and $\pi_{117}(L_3V(2))$ are all trivial.}

\p
In the $E(3)$-based Adams spectral sequence \kko{EnASS}, the $E_2$-term $E_2^{*,*}(V(2))$ for the homotopy groups in the lemma are given by
$$
H^{5,120}K(3)_*,\ak  H^{4,120}K(3)_*\qand  H^{3,120}K(3)_*.
$$
Then, $\rank H^{s,t}K(3)_*\le \rank (K(3)_*\ox H^*L(3,3))$ for the module $H^*L(3,3)$  determined in \cite[(3.8) Th.]{r}:
\AL{
H^3L(3,3)&=A^2\ze_3\o+ A^3,\\
H^4L(3,3)&=A^3\ze_3\o+ A^4 \qand\\
H^5L(3,3)&=A^4\ze_3\o+ (A^3\ze_3)^*\o+ (A^4)^*
}
for
\AL{
A^2&=\Z/5\{g_i, k_i, \b1i\},\\
A^3&=\Z/5\{g_i\h1{i+1}, \el1i, \el2i, \el3i, \el4i, \el5i\}\qand\\
A^4&=\Z/5\{m_{i,j}, m'_i\}.
}
Here,
\AR {ll}{
\el1i=\h1i\h2i\h3i, & \el2i=\h1i\h2i\h2{i+2},\\ \el3i=\h1i\h2i\h2{i+1}+\h1i\h1{i+1}\h3i,&
\el4i=\h1i\h2{i+2}\h3{i+1},\\
 \el5i=\sum_i(\h1i\h2{i+1}-\h1{i+1}\h2{i+2})\h3i
}
and 
\AL{
m_{i,j}&=\h1ik_i\h3j=g_i\h1{i+1}\h3j\qand\\
m'_i&= \h1{i+2}\h1i\h2i(\h3i+\h3{i+1})\pm \h1i\h20\h21\h22.
}
These elements have the degrees (modulo $q_3=248=|v_3|$) as follows:
\AR{|c|c|c|c|c|c|c|c|c|}{
\hline
56& 88& 40& 32&192&200&160&-32&8\\
\hline
g_0&k_0&\b10&g_1&k_1&\b11&g_2&k_2&\b12\\
\hline
}
\AR{|c|c|c|c|c|c||c|c|}{
\hline
96&56&16&48&-32&0&96&8\\
\hline
g_0\h11&\el10&\el20&\el30&\el40&\el50&m_{0,j}&m_0'\\
\hline\hline
232&32&80&240&88&0&232&40\\
\hline
g_1\h12&\el11&\el21&\el31&\el41&\el51&m_{1,j}&m_1'\\
\hline\hline
168&160&152&-40&192&0&168&200\\
\hline
g_2\h10&\el12&\el22&\el32&\el42&\el52&m_{2,j}&m_2'\\
\hline
}
The dual degrees of these elements (elements in $(A^3\ze_3)^*\mid (A^4)^*$) are
\AR {cccccc|cc}{
152,& 192,& 232,& 200,& 32,& 0&152,&240;\\
16,& 216,& 168,&8,& 160,& 0&16,&208;\\
80,& 88,& 96,& 40,& 56,& 0&80,&48.
}
Thus, 
there is no element with degree 120. \q

\subsection{The case for $(p,n)=(7,4)$:}

The dimensions of cells in $V(3)$ are:
\AL{
\Zt(3)=\{0, 1, 13,14, 97,98,110,111, 685, 686, 698,699,782,783,795,796\}.
}
In \cite[Th.~4.4]{toda}, Toda showed the existence of the pairing $V(2\frac14)\sm V(3)\to V(3)$, which follows from the fact 
$$
\pi_{i-1}(V(3))=0\mbx{for $i=s+a>1$ with $s\in \Zt(3)\jo{686}$ and $a\in \Zt(3)$.} \lnr{V30}
$$
Here, $\Zt(3)\jo t=\{s\in \Zt(3)\mid s\le t\}$.

Let $W$ be a spectrum sitting in the cofiber sequence
$$
\Si^{738}S^0\xar{\be_1^9} S^0\xar{i_W}W\xar{j_W} \Si^{739}S^0, \lnr{cofW}
$$
in which $\be_1\in \pi_{82}(S^0)$ is the well known generator.
Then, $W$ is a ring spectrum by \cite[Cor.~2.6]{oka}.

We show the existence of the pairing $\ph'\cln V(3)\sm V(3)\to V(3)\sm W$ in Proposition \ref{V3W}, and $\be_1^8=0\cln V(3)\to L_4V(3)$ in Lemma \ref{be190} below.
The lemma implies the decomposition
 $L_4V(3)\sm W=L_4V(3)\vee \Si^{739}L_4V(3)$.
Therefore, we obtain the composite
$$
V(3)\sm V(3)\xar{\ph'} V(3)\sm W \to L_4V(3)
$$
for $\ph'$ in \kko{ph'},
which yields the desired ring structure on $L_4V(3)$. 

We now prove Proposition \ref{V3W} and Lemma \ref{be190}.

\lem{V3W0}{The homotopy groups
$
\pi_{i-1}(V(3)\sm W) $ are trivial for $i=s+a$ with $s,a\in \Zt(3)\setminus \Zt(3)\jo{686}$.
}

\p
Since we have an exact sequence
\ARr c{
\pi_{i-739}(V(3))\xar{\be_1^9} 
\pi_{i-1}(V(3))\xar{(i_W)_*}\pi_{i-1}(V(3)\sm W)\\
\hspace{2in} \xar[(j_W)_*]{} \pi_{i-740}(V(3))\xar[\be_1^9]{} \pi_{i-2}(V(3)), 
}{exactW}
we study the homotopy groups $\pi_{i}(V(3))$ under \kko{piV}.
Table \ref{table} gives rise to the following table of 
$\LR{H^{*}U(L)\ox\Z/7\{\b11,\b20\}}^{s,t}$ with $t-s\le (2p^2+4p+6)q+7=1591$:
{\small
\AR{|c||c|c|c|c|c|c|c|c|c|}{
\hline
\!\!\!\!\mbox{\tiny $\|x\|$}\!\!\!\!&0&11&83(1)&106(24)&178(14)&189(25)&586(12)&587(13)\\
\hline
\!\!x\!\!&1&h_0&h_1&g_0&k_0&k_0h_0&\b11&h_2\\
\hline\hline
\!\!\mbox{\tiny $\|x\|$}\!\!&597(23)&598(24)&669(13)&670(14)&681(25)&692(36)&753(15)&754(16)\\
\hline
\!\!x\!\!&\b11h_0&h_2h_0&\b11h_1&\b20&\b20h_0&\b11g_0&\b20h_1&g_1\\
\hline\hline
\!\!\mbox{\tiny $\|x\|$}\!\!&764(26)&775(37)&776(38)&789(51)&848(28)&849(29)&859(39)&872(52)\\
\hline
\!\!x\!\!&\b11k_0&\b11k_0h_0&\b20g_0&l_1&\b20k_0&l_2&\!\!\b20k_0h_0\!\!&l_1h_1\\
\hline\hline
\!\!\mbox{\tiny $\|x\|$}\!\!&1172(24)&1173(25)&1183(35)&1184(36)&1255(25)&\!\!1256(26)\!\!&\!\!1257(27)\!\!&\!\!1258(28)\!\!\\
\hline
\!\!x\!\!&\b11^2&\b11h_2&\b11^2h_0&\b11h_2h_0&\b11^2h_1&\b11\b20&\b20h_2&k_1\\
\hline
\hline
\!\!\mbox{\tiny $\|x\|$}\!\!&1267(37)&1268(38)&1278(48)&1281(51)&1339(27)&\!\!1340(28)\!\!&\!\!1340(28)\!\!&\!\!1341(29)\!\!\\
\hline
\!\!x\!\!&\!\!\b11\b20h_0\!\!&\b20h_2h_0&\b11^2g_0&l_3&\!\!\b11\b20h_1\!\!&\b20^2&\b11g_1&k_1h_1\\
\hline
\hline
\!\!\mbox{\tiny $\|x\|$}\!\!&1350(38)&1351(39)&1361(49)&1362(50)&1375(63)&\!\!1376(64)\!\!&\!\!1423(29)\!\!&\!\!1424(30)\!\!\\
\hline
\!\!x\!\!&\b11^2k_0&\b20^2h_0&\b11^2k_0h_0&\b11\b20g_0&\b11l_1&l_1h_2&\b20^2h_1&\b20g_1\\
\hline
\hline
\!\!\mbox{\tiny $\|x\|$}\!\!&1434(40)&1435(41)&1445(51)&1446(52)&1458(64)&\!\!1459(65)\!\!&\!\!1518(42)\!\!&\!\!1519(43)\!\!\\
\hline
\!\!x\!\!&\!\!\b20\b11k_0\!\!&\b11l_2&\!\!\b20\b11k_0h_0\!\!&\b20^2g_0&\b11l_1h_1&\b20l_1&\b20^2k_0&\b20l_2\\
\hline\hline
\!\!\mbox{\tiny $\|x\|$}\!\!&\!\!1529(53)\!\!&1532(56)&1542(66)&1543(67)\\
\cline{1-5}
\!\!x\!\!&\b20^2k_0h_0&m_1&\b20l_1h_1&m_1h_0\\
\cline{1-5}
\multicolumn{9}{c}{\mbox{\tiny }}\\
\multicolumn{9}{c}{\mbox{\normalsize Table \nr \label{table7}
}}
}
}%
Here, in the rows $\|x\|$, the numbers $w(u)$ denote the total degrees of the elements $x$ under them:
$$
w=\|x\|=t-s ,\ak u\cg w \mod 82\qand 0\le u<82,
$$
in which $82=|\b10|-2=\|\be_1\|$.
In order to find a generator $\LR{P(\b kl)\ox H^{*,*}(U(L))}^{s,t}$ with $i=t-s$,
it suffices to find $w(u)$ in the $\|x\|$ rows in Table \ref{table7} such that 
\begin{itemize}
\item
$i\cg u$ mod $82$ and 
\item $w\le i$. 
\end{itemize}
If we find such $w(u)$, then $\LR{P(\b kl)\ox H^{*,*}(U(L))}^{s,t}$ contains an element of the form
$$
x\b10^c
$$
for the integer $c=\frac{i-w}{82}$,
where $x$ is the element under $w(u)$ in the table.

Furthermore, we notice the  relation 
$$
\|x\|=\|y\|+1 \mbx{if $d_r(x)=y$} 
$$ of the total degrees in all of the May spectral sequences $P(\b kl)\ox H^{*,*}(U(L))\Rightarrow H^*(V(L))$ and
$H^*(V(L))\Rightarrow H^*\cP\cong E_2^*(V(3))$ (\cf \cite{toda}), and the Adams-Novikov spectral sequence $E_2^*(V(3))\Rightarrow \pi_*(V(3))$.
This implies 

\medskip

\noindent
(\nr) \label{perm} {\it
Consider an element $x$ with total degree $w(u)$ in Table \ref{table7} and
suppose that $x$ yields a permanent cycle in the $E_2$-term $E_2^*(V(3))$.
Then, if $x$ does not survive to the homotopy group $\pi_*(V(3))$, then it is the image of an element of total degree $w+1(u+1)$ under some differential of the above spectral sequences. 
In particular, $x$ yields an essential homotopy element of
$\pi_w(V(3))$ if there is no element with degree $w+1(u+1)$.
}

\medskip

For the integers in $\Zt(3)\setminus \Zt(3)\jo{686}$, we notice the following:
\ARr{|c||c|c|c|c|c|c|}{
\hline
& 698&699&782&783&795&796\\
\hline
\mx{\mbox{mod}\\
(82)}
&42&43&44&45&57&58\\
\hline
}{n3}

We first consider $\pi_{s+a-740}(V(3))$ for $s, a\in \Zt(3)\setminus \Zt(3)\jo{686}$.
Then, the integers $s+a-740$ are: 
\AR l{
656(0), \ak 657(1),\ak 740(2), \ak 741(3), \ak 753(15), \ak 754(16);\\
 658(2),\ak 741(3), \ak 742(4), \ak 754(16), \ak 755(17);\\
824(4), \ak 825(5), \ak 837(17), \ak 838(18); \ak 826(6), \ak 838(18), \ak 839(19);\\
850(30), \ak 851(31); \ak 852(32),
}
and we find elements
$$\b10^8\ (656(0)),\ak 
h_1\b10^7\  (657(1)), \ak \b20h_1 \ (753(15)) \qand g_1\ (754(16)) \lnr{genW}
$$
from Table \ref{table7}.

Next consider similarly $\pi_{s+a-1}(V(3))$ for $s, a\in \Zt(3)\setminus \Zt(3)\jo{686}$.
Then, the integers $s+a-1$ are: 
\AR l{
1395(1),\ak 1396(2), \ak 1479(3), \ak 1480(4),\ak 1492(16),\ak 1493(17);\\
1397(3), \ak 1480(4), \ak 1481(5),\ak 1493(17),\ak 1494(18);\\
1563(5), \ak1564(6), \ak 1576(18), \ak 1577(19);\\
1565(7), \ak 1577(19), \ak 1578(20);\ak 1589(31),\ak 1590(32);\ak  1591(33).
}
Thus, we find
$$
h_1\b10^{16}\ (1395(1))\qand g_1\b10^9\ (1492(16)). \lnr{obst}
$$
(It is stated in \cite[Th.~4.4]{toda} that $h_1\b10^{2p+3}$ is an obstruction, but it is $h_1\b10^{2p+2}$ as shown above.)

Let $\llb x$ denote the homotopy element detected by $x$. For example, $\llb{\b10}=\io_3\be_1\in \pi_{82}(V(3))$. 
Hereafter, $\io_3$ denotes the inclusion
$$
\io_3=i_{I_4}=i_3i_2i_1i\cln S^0\xar i V(0)\xar {i_1}V(1)\xar {i_2}V(2)\xar {i_3}V(3) \lnr{io3}
$$
to the bottom cell.
We notice that\\[-1mm]

\noindent
(\nr) \label{non0} {\it
the elements $\b10^8$, $\b10^{17}$, $h_1\b10^7$ and $h_1\b10^{16}$ detect essential homotopy elements} \\[-1mm]

\noindent
by \kko{perm},
since $h_1$ and $\b10$ detect the generators $\be_1'\in\pi_{83}(V(0))$ and $\be_1\in\pi_{82}(S^0)$, respectively.
Therefore, $\be_1^9(\llb {\b10^8})=\llb {\b10^{17}}\ne 0$ and $\be_1^9(\llb {h_1\b10^7})=\llb {h_1\b10^{16}}\ne 0$ by \kko{non0}, and so in \kko{exactW},
$$
\llb{\b10^8},\ \llb{h_1\b10^7}\not\in\Ker \be_1^9\qand \llb{h_1\b10^{16}}\in \im \be_1^9. \lnr{bh}
$$

We next consider the elements $\b20h_1$, $g_1$ and $g_1\b10^9$ in \kko{genW} and \kko{obst}.
Note that $g_1=\lrk{h_1,h_1,h_2}\in E_2^{2,*}(V(3))$.
By virtue of \kko{perm} together with Table \ref{table7}, $d_r(g_1)\in E_r^{r+2, 755+r}(V(3))$ must be of the form $\b10^i\b20h_1$, but this is not the case by degree reason.
Therefore, $g_1$ is a permanent cycle and $g_1\b10^9$ detects an essential homotopy element by \kko{perm} with Table \ref{table7}.
Furthermore, this also implies that $\b10^9\b20h_1$ is not a target of the differential by \kko{perm}, if $\b20h_1$ is a permanent cycle. 
Thus,
$$
\llb{\b20h_1},\ \llb{g_1}\not\in \Ker \be_1^9\qand \llb{h_1\b10^{16}}\in \im \be_1^9. \lnr{g1}
$$
Now the lemma follows from \kko{exactW}, \kko{genW}, \kko{obst}, \kko{bh} and \kko{g1}.
\q

\prop{V3W}{
$V(3)\sm W$ is a ring spectrum.
}

\p
We first show the existence of a pairing $\ph'\cln V(3)\sm V(3)\to V(3)\sm W$ such that $ \ph'(\io_3\sm V(3))=V(3)\sm i_W$ for $\io_3$ in \kko{io3} by attaching cells as Toda did. 
By \kko{V30}, we have an extension $V(3)\jo{697}\sm V(3)\to V(3)$ of the identity $V(3)\to V(3)$.
Thus, we have a composite $V(3)\jo{697}\sm V(3)\to V(3)\xar{V(3)\sm i_W}V(3)\sm W$.
Lemma \ref{V3W0} certifies the existence of an extension  
$$
\ph'\cln V(3)\sm V(3)\to V(3)\sm W\lnr{ph'}
$$ of the composite. 

Now the multiplication of $V(3)\sm W$ is given by
\AL{
\LR{V(3)\sm W}\sm \LR{V(3)\sm W}\xar[\simeq]{1\sm T\sm 1}&\ V(3)\sm V(3)\sm W\sm W\xar{\ph'\sm \mu_W}
V(3)\sm W\sm W\\
&\xar{1\sm \mu_W} V(3)\sm W.
}
Here, $T$ denotes the switching map and $\mu_W$ denotes the multiplication of the ring spectrum $W$.
\q

Assuming Corollary \ref{be140}, which is shown independently, we see the following:

\lem{be190}{
 $\be_1^8\sm L_4V(3)=0\in [V(3),L_4V(3)]_{656}$.
}

\p
Consider the cofiber sequence
$$
\Si^{684}V(2)\xar{\ga} V(2)\xar{i_3}V(3)\xar{j_3}\Si^{685}V(2).
$$
Since we have a pairing
$\ph_{23}\cln V(2)\sm V(3)\to V(3)$ such that $i_3=\ph_{23}(V(2)\sm \io_3)$
for the inclusion $\io_3\cln S^0\to V(3)$ in \kko{io3},
\AL{
\be_1^s\sm i_3&=\be_1^s\sm (\ph_{23}(V(2)\sm \io_3))=\ph_{23}(V(2)\sm \io_3\be_1^s)\cln \Si^{82s}V(2)\to V(3).
}
By Corollary \ref{be140}, $\io_3\be_1^4=0$, and so $\be_1^4\sm i_3=0$.
Consider a commutative diagram
{\small
\Xy{8}{
[V(3), L_4V(3)]_*\ar[d]_-{i_3^*}\ar[r]^-{ad}_-\cong &[V(3)\sm DV(3),L_4S^0]_*\ar[d]^-{i_3^*}\ar[r]^-{ad}_-\cong& [DV(3), L_4DV(3)]_*=[V(3), L_4V(3)]_{*}\ar@<15mm>[d]^-{(j_3)_*}\ar@<-13mm>[d]^-{D(i_3)_*}\\
[V(2), L_4V(3)]_*\ar[r]^-{ad}_-\cong &[V(2)\sm DV(3),L_4S^0]_*\ar[r]^-{ad}_-\cong& [DV(3), L_4DV(2)]_*=[V(3), L_4V(2)]_{*},
}
}%
in which $ad$ denotes the adjunction.
This together with the relation $\be_1^4\sm i_3=0$ gives rise to $\be_1^4\sm j_3=0$.

Therefore, we have elements $\xi_4\in [V(2),L_4V(3)]_{1013}$ and $\xi_4^*\in [V(2),L_4V(3)]_{328}$ such that
 $$
\be_1^4\sm L_4V(3)=\xi_4j_3\qand =i_3\xi_4^* : 
$$
\xyR 9{
&\Si^{656}L_4V(3)\ar[d]^-{\be_1^4\sm 1}\ar@{..>}[dl]_-{\xi_4^*}\\
\Si^{328}L_4V(2)\ar[r]^-{i_3}&\Si^{328}L_4V(3) \ar[r]^-{j_3}\ar[d]_-{\be_1^4\sm 1}&\Si^{1013}L_4V(2)\ar@{..>}[dl]^-{\xi_4}\\
&L_4V(3)
}
Thus,
$$
\be_1^8\sm V(3)=\xi_4j_3i_3\xi^*_4=0
$$
as desired. \q

\section{Proof of Theorem \ref{small}}

Put
$$
H^{s,t}M=\e_{E(n)_*(E(n))}^{s,t}(E(n)_*,M)
$$
for an $E(n)_*(E(n))$-comodule $M$.
Then, we have the Miller-Ravenel change of rings theorem 
$$
\e_{BP_*(BP)}^*(BP_*, v_n^{-1}BP_*/I_n)\cong
H^*(E(n)_*/I_n)=H^*(K(n)_*)
$$ for the invariant ideal $I_n=(p,v_1,\dots, v_{n-1})$ (\cf \cite[Th.~3.1]{hs}).
Here, $K(n)_*=E(n)_*/I_n=\Z/p[v_n,v_n^{-1}]$.
Ravenel introduced in  \cite[\S1] r the exterior complex
$$
C(n)=E(h_{i,j}\colon 1\le i\le n,~j\in\mathbb Z/n)
$$
with differential given by
$$
d(h_{i,j})=\sum_{\ell=1}^{i-1}h_{\ell,j}h_{i-\ell,\ell+j}. \lnr{diff}
$$
Here,
the bidegree of the generator $h_{i,j}$ is $(1,p^jq_i)\in \Z\times (\Z/q_n)$ for $q_i$ in \kko{qn}. 
From the results in  \cite[\S1] r,
 we deduce the following:

\lem{coh} {Let $J$ be an invariant ideal $(p,v_1^{e_1},\dots,v_{n-1}^{e_{n-1}})$ of $E(n)_*$. Then, $$
\rank H^{s,t}E(n)_*/J\le \rank \LR{E(n)_*/J\ox H^*(C(n),d)}^{s,t}
$$
as $\Z/p$-modules.
In particular, $\rank H^{s,t}K(n)_*\le \rank \LR{K(n)_*\ox H^*(C(n),d)}^{s,t}$.
}

Let  $M^{*,*}$ denote a basis of the $\mathbb Z/p$-vector space $C(n)$ consisting of monomials.
In the following, we use the word ``monomial" for an element of $M^{*,*}$.
Note that  $M^{n^2,*}$ consists of only one element of degree $0$.
We denote the element of  $M^{n^2,0}$ by $g$.
For $x\in M^{s,*}$, we define $x^*\in M^{n^2-s,*}$ by
$$
xx^*=\pm g,\lnr{star}
$$
and obtain an  isomorphism
$$
C(n)^{q+1,tq}\xrightarrow[\cong]{(-)^*} C(n)^{n^2-q-1,-tq}
\lnr{q+1,0}
$$
given
by $(-)^*(x)=x^*$ for $x\in M^{q+1,tq}$. 

\prop{Hn2}{Suppose $n<p-1$ and let $J_k=I_{n-1}+(v_{n-1}^k)$. Then,
$H^{n^2,q}E(n)_*/J_k=0$ for each $1\le k\le \sum_{i=2}^{n-1}p^i$.
}

\p
Since $|g|=0=|v_n|$, we will find a positive integer $a$ such that  $|v_{n-1}^a|=q$.
This gives an equation
$$
a(p^{n-1}-1)=p-1+b(p^{n}-1)\in \Z
$$
for an integer $b$.
It follows that $a=bp+(b+1)/\LR{\sum_{i=0}^{n-2}p^i}$ and so $b+1=u\sum_{i=0}^{n-2}p^i$ for some $u\ge 1$.
Therefore,
$$
a=u\sum_{i=0}^{n-1}p^i-p
$$
for $u\ge 1$.
Thus, if $k<\sum_{i=0}^{n-1}p^i-p$, there is no generator in $H^{n^2,q}E(n)_*/J_k$ by Lemma \ref{coh}.
\q

\noindent
{\it Remark.} More careful computation using \cite[(5.18)]{mrw} makes $k$ in the proposition greater. \\

We note that 
$$
E_2^{q+1,q}(V(n-1))=H^{q+1,q}E(n)_*/I_n=H^{q+1,q}K(n)_*.\lnr{E2HK}
$$

\cor{pn53}{
For $(p,n)=(5,3)$, $H^{9,8}E(3)_*/J_k=0$ for $k\le 25$.
In particular, Theorem \ref{small} holds for $(p,n)=(5,3)$. 
}

Now turn to the case $(p,n)=(7,4)$.

\lem{C(n)}{ Let $(p,n)=(7,4)$.
\begin{enumerate}
\item
$C(4)^{3, -12}$ is generated by the elements
\AL{
&\h31\h4i\h4j \mbx{for $0\le i<j\le 3$,}\\
 &\h11\h22\h4i\qand \h13\h21\h4i \mbx{for  $0\le i\le 3$,} \\
&\h11\h12\h13, \ak \h11\h31\h32, \ak \h12\h31\h33,\ak \h13\h30\h31, \\
&\h20\h22\h31, \ak \h21\h22\h33\qand \h21\h23\h31.
}

\item
$C(4)^{8,336}=0$.
\end{enumerate}
}

\p
Consider the subalgebra
$$
\O C(4)=E(\h ij:1\le i\le 3,\  j\in \Z/4) \lnr{subal}
$$
of $C(4)$.
The degree $|x|\in \Z/q_4=\Z/4800$ of a monomial $x\in \O C(4)$ is expressed
as
$$
|x|=12\times \sum_{i=0}^37^ia_{4-i} \mbx{with $0\le a_i\le 6$,}
$$  
which we write
$$
|x|=(a_1a_2a_3a_4).
$$
We further assume that the integers $a_i$ satisfy $(a_1a_2a_3a_4)\le (1111)$ under the lexicographic order. 
We also use a similar notation
$$
(a_1a_2a_3a_4)_\N=12\times \sum_{i=0}^37^ia_{4-i}\mbx{for $a_i\ge 0$.}
$$
Note that 
$$
|x|\cg ((a_1+k)(a_2+k)(a_3+k)(a_4+k))_\N \mod 4800 \mbx{for an integer $k\ge 0$.}
$$
For a monomial $x\in \O C(4)$, we also introduce notations 
$$
 [x]=a_1+a_2+a_3+a_4 \qand  (x)_i=a_{4-i}
$$ if $|x|=(a_1a_2a_3a_4)_\N$. 

For the algebraic generators $\h ij$ of $\O C(4)$, we have\\

\noindent
(\nr) \label{alg-gen}
\begin{minipage}{4.9in}
\begin{enumerate} \renewcommand{\labelenumi}{\alph{enumi})}
\item \ak $\h 11$ $(0010)$, $\h 20$ $(0011)$, $\h 21$ $(0110)$,
 $\h 30$ $(0111)$, $\h 31$~$(1110)$, \\ \phantom{$\h 31$~$(1110)$,} $\h 33$ $(1011)$; 
\item \ak $\h 10$ $(0001)$, $\h 12$ $(0100)$, $\h 13$ $(1000)$,
 $\h 22$ $(1100)$, $\h 23$~$(1001)$, \\ \phantom{$\h 31$~$(1110)$,}  $\h 32$ $(1101)$.
\end{enumerate}
\end{minipage}

\medskip

The generators $\h ij$ in a) and b) satisfy $(\h ij)_1=1$ and $(\h ij)_1=0$, respectively.

\begin{enumerate}
\item We will find monomials $x\in \O C(4)^{s,-12}$ with $s\le 3$.
Then 
$$
|x|=(kkk(k-1))_\N\mbx{for an integer $k\ge 1$.}
$$
If $k=1$, that is, $(1110)_\N$, then from \kko{alg-gen}, we find
$$
\h31\ (s=1),\ak \h21\h13, \ak \h11\h22\ (s=2), \qand \h11\h12\h13\ (s=3).
$$
Turn to the case $k=2$ ($|x|=(2221)_\N$ and $[x]=7$).
If $\h1i$ is a factor, say $x=\h1i y$, then $[y]=6$, and so $y=\h3{j}\h3{j'}$:
$$
\h13\h30\h31,\ak \h12\h31\h33,\ak \h11\h31\h32.
$$
If $\h2i$ is a factor, say $x=\h2i y$, then $[y]=5$, and so  $\h2j\h3{j'}$:
$$
\h20\h31\h22, \ak \h21\h33\h22,\ak \h22\h33\h21,\ak \h23\h31\h21.
$$
For $k\ge 3$, then $[x]=4k-1\ge 11$, while $[\h ij\h {i'}{j'}\h {i''}{j''}]\le 9$.
Therefore, no element satisfies this. 

Now the first statement of the  lemma follows from the relation $$C(4)^{3,-12}=\O C(4)^{3,-12}\o+\Op_{i}\h4i\O C(4)^{2,-12}\o+ \Op_{i,j}\h4i\h4j\O C(4)^{1,-12}.$$

\item
We will find monomials $x\in \O C(4)^{s,336}$ with $4\le s\le 8$.
In this case, it may have a carry-over:
$|x|=(0037), \dots$.
Since there are at most six algebraic generators $\h ij$ of $\O C(4)$ with $(\h ij)_k=1$ for each $k$,
none of the carry-over case occurs. 
Thus, we consider an element $x\in \O C(4)$ with 
$$
|x|\cg (kk(k+4)k)_\N \mod 4800 \mbx{for an integer $0\le k\le 2$.}
$$
This implies that $x$ has a factor of the form $\h {i_1}{j_1}\h {i_2}{j_2}\h {i_3}{j_3} \h {i_4}{j_4}$ for $\h{i_l}{j_l}$
in  \kko{alg-gen} a).
Therefore, $[x]\ge 8$, and so $k\ge 1$.
Therefore, $k=1,2$.

If $k=1$, then $[x]=8$, and  $[\h 11\h 20\h21\h3i]=8$.
Then, 
\AL{
|\h 11\h 20\h21\h3i|=&(0242)_\N \ (i=0), (1241)_\N\ (i=1),\\
&(1232)_\N\ (i=2),(1142)_\N\ (i=3).
}
Thus, these yield no solution.

If $k=2$, then $[x]=12$ and $(x)_1=6$.
Since $(x)_1=6$, $x$ has all elements in \kko{alg-gen} a) as a factor.
Let $h$ be the product of the six elements in \kko{alg-gen} a).
Then, $[h]=14$. This contradicts to $[x]\ge [h]$.

Thus, we have the second from 
\AL{
C(4)^{8,336}=&\O C(4)^{8,336}\o+\Op_{i}\h4i\O C(4)^{7,336}\o+\Op_{i,j}\h4i\h4j\O C(4)^{6,336}\\
&\o+\Op_{i,j,k}\h4i\h4j\h4k\O C(4)^{5,336}\o+\h40\h41\h42\h43\O C(4)^{4,336}.
} 
\end{enumerate}

\vspace{-.2in}

\q

From Lemma \ref{C(n)} 2), we deduce

\cor{be140}{$\pi_{328}(L_4V(3))=0$, and in particular, $\io_3\be_1^4=0\in \pi_{328}(L_4V(3))$ for the inclusion $\io_3$ in \kko{io3}. }

\p
By the spectral sequences,
$$
\rank \pi_{328}(L_4V(3))\le \rank E_2^{8,336}(V(3))\le \rank(K(4)_*\ox H^*(C(4),d))^{8,336}.
$$
Furthermore, the right hand side is also not greater than 
rank $(K(4)_*\ox C(4))^{8,336}$, which is $0$ by Lemma \ref{C(n)}.
Therefore, we obtain $\pi_{328}(L_4V(3))=0$. 
\q

\noindent
{\it Proof of Theorem \ref{small} for $(p,n)=(7,4)$.}
By Lemma \ref{C(n)} together with \kko{q+1,0},
$C(4)^{13,12}$ is generated by the set $M_1\cup M_2\cup M_3(\sus M^{*,*})$ of  monomials given by
\AL{
M_1&=\{g(i,j)=(\h31\h4i\h4j)^*\mid 0\le i<j\le 3\}, \\
M_2&=\{g_1(i)= (\h11\h22\h4i)^*,\ g_2(i)= (\h13\h21\h4i)^*\mid 0\le i\le 3\}\qand \\
M_3&=\{g_1= (\h11\h12\h13)^*,\ g_2= (\h11\h31\h32)^*, \ g_3= (\h12\h31\h33)^*,\\
&\qquad \ak  g_4=(\h13\h30\h31)^*, \
 g_5= (\h20\h22\h31)^*, \ g_6= (\h21\h22\h33)^*\\
&\hspace{1in}  g_7= (\h21\h23\h31)^*\}.
}
On the differential $d$,  we notice that 
\naka{
$d(x^*)\doteq \sum y^*$ for monomials $x$ and $y$ if and only if $d(y)\doteq x+\dots$.
}
Here, $\doteq$ denotes equality up to sign.
In particular, $d((\h ij\h kl)^*)=0$ if $i+k\ge 5$.
Under these facts with \kko{diff}, we obtain
{\small
\AL{
d((\h11\h22\h4i\h4j)^*)&\doteq (\h31\h4i\h4j)^*\doteq g(i,j),\\
d( (\h11\h22\h4i)^*)&\doteq (\h31\h4i)^*\qquad (d(g_1(i))\doteq (\h31\h4i)^*),\\
d((\h11\h12\h13\h4i)^*)&\doteq (\h13\h21\h4i)^*+(\h11\h22\h4i)^*\doteq g_1(i)+g_2(i),\\
d( (\h11\h12\h13)^*)&\doteq  (\h21\h13)^*+(\h11\h22)^*  \
\ak (d(g_1)\doteq  (\h21\h13)^*+(\h11\h22)^*),\\
d((\h10\h13\h21\h31)^*)&\doteq (\h23\h21\h31)^*+(\h13\h30\h31)^*+(\h13\h21\h40)^*+(\h13\h21\h41)^*\\
&\doteq  g_7+g_4+g_2(0)+g_2(1),\\
d((\h10\h11\h22\h31)^*)&\doteq (\h20\h22\h31)^*+(\h32\h11\h31)^* +(\h11\h22\h40)^*+(\h11\h22\h41)^*\\
&\doteq  g_5+g_2+g_1(0)+g_1(1),\\
 d((\h11\h12\h22\h33)^*)&\doteq (\h21\h22\h33)^*+(\h31\h12\h33)^*+(\h11\h22\h42)^*+(\h11\h22\h43)^*\\
&\doteq  g_6+g_3+g_1(2)+g_1(3),\\
d((\h11\h21\h22\h23)^*)&\doteq (\h21\h23\h31)^*+(\h21\h22\h33)^* +(\h11\h22\h41)^*+(\h11\h22\h43)^*\\
&\doteq  g_7+g_6+g_1(1)+g_1(3),\\
d((\h13\h20\h21\h22)^*)&\doteq (\h21\h22\h33)^*+(\h20\h22\h31)^*+(\h13\h21\h40)^*+(\h13\h21\h42)^*\\
&\doteq  g_6+g_5+g_2(0)+g_2(2)\qand\\
d((\h11\h12\h23\h31)^*)&\doteq (\h11\h31\h32)^*+(\h21\h23\h31)^*+(\h12\h33\h31)^*\\
&\doteq  g_2+g_7+g_3.
}%
These give rise to an exact sequence
$$
0\to A\xar{d} C(4)^{13,12}\xar{d} B\to 0
$$
for the submodules $A\sus C(4)^{12,12}$ and $B\sus C(4)^{14,12}$ given by 
\AL{
A&=\Z/7\{(\h11\h22\h4i\h4j)^*,(\h11\h12\h13\h4k)^*,(\h10\h13\h21\h31)^*,\\
&\qquad  \qquad  (\h10\h11\h22\h31)^*, (\h11\h12\h22\h33)^*, (\h11\h21\h22\h23)^*,\\
&\qquad   \quad  (\h11\h12\h23\h31)^* ,(\h13\h20\h21\h22)^*
\mid 0\le i<j\le 3, 0\le k\le 3\}\qand \\
B&=\Z/7\{(\h31\h4i)^*,  (\h21\h13)^*+(\h11\h22)^*\mid 0\le i\le 3\}.
}
Indeed, the images of the generators of $A$ under $d$ are linearly independent, and $\rank A=16$, $\rank B=5$ and $\rank C(4)^{13,12}=21$.
This implies $H^{13,12}C(4)=0$ and then the theorem by Lemma \ref{coh} with \kko{E2HK}.
\q

\vspace{-.5cm}

\bibliographystyle{amsplain}

\end{document}